\documentclass[11pt, a4paper]{article}

\usepackage{amssymb}
\usepackage{amsmath}
\usepackage{amscd}
\usepackage{theorem}
\usepackage{mathrsfs}
\usepackage{enumerate}
\usepackage{graphicx}
\usepackage[T1]{fontenc}
\usepackage{ae}
\usepackage{pictexwd,dcpic}

\theoremstyle{change}

\newtheorem{thm}{Theorem.}
\newtheorem{lem}[thm]{Lemma.}

\newtheorem{cor}[thm]{Corollary.}

\renewcommand{\phi}{\varphi}

\newcommand{\proof}{\textbf{Proof. }}
\newcommand{\proofend}{\hfill $\Box$}

\begin{document}

\centerline{\textbf{\large Multivariable Bergman shifts and Wold decompositions}}

\vspace{1cm}

\centerline{J\"org Eschmeier and Sebastian Langend\"orfer }

\vspace{1cm}

\begin{center}
\parbox{12cm}{\small 
Let $H_m(\mathbb B)$ be the analytic functional Hilbert space on the unit ball $\mathbb B \subset \mathbb C^n$
with reproducing kernel $K_m(z,w) = (1 - \langle z,w \rangle)^{-m}$. 
Using algebraic operator identities we characterize those commuting row contractions $T \in L(H)^n$ on a Hilbert space $H$
that decompose into the direct sum of a spherical coisometry and copies of the
multiplication tuple $M_z \in L(H_m(\mathbb B))^n$. For $m=1$,
this leads to a Wold decomposition for partially isometric commuting row contractions that are
regular at $z = 0$. For $m = 1 = n$, the results reduce to the classical
Wold decomposition of isometries. We thus extend corresponding one-variable results 
of Giselsson and Olofsson \cite{GO} to the case of the unit ball.

\vspace{0.5cm}

}

\end{center}

\centerline{\textbf{\S1 \, Introduction}}

\vspace{.5cm}

By the classical {\it Wold decomposition} theorem each isometry $T\in L(H)$ on a Hilbert space $H$ 
is a direct sum $T = T_0\oplus T_1 \in L(H_0 \oplus H_1)$ of
a unitary operator $T_0 \in L(H_0)$ and an operator $T_1 \in L(H_1)$ which is unitarily 
equivalent to a Hardy space shift $M_z \in L(H^2(\mathbb D,\mathcal D))$. Isometries are characterized by the 
the operator identity $T^\ast T = 1_H$, and the Hardy space is the analytic functional Hilbert space on the unit
disc with reproducing kernel $K(z,w) = (1 - z\overline{w})^{-1}$. Our aim is to prove corresponding 
decomposition theorems for commuting tuples $T \in L(H)^n$ of Hilbert space 
operators which satisfy higher order operator identities related to the reproducing kernel 
$K_m(z,w) = (1-\langle z,w\rangle)^{-m}$ on the unit ball.

An operator $T \in L(H)$ on a Hilbert space $H$ is unitarily equivalent to a Hardy space shift $M_z \in L(H^2 (\mathbb D,\mathcal D))$ 
if and only if it is an isometry 
which is {\it pure} in the sense that $\bigcap^\infty_{k=0}T^k H$ $= \lbrace 0 \rbrace$. We replace the Hardy space $H^2(\mathbb D)$ on the unit disc by the analytic functional Hilbert spaces 
$H_m(\mathbb B)$ on the open unit ball $\mathbb B \subset \mathbb C^n$ defined by the reproducing kernels 
$K_m(z,w)=(1-\langle z,w\rangle)^{-m}$, where $m \geq 1$ is a positive integer. It is well known 
that the multiplication tuple $M_z=(M_{z_1},\ldots,M_{z_n})\in L(H_m(\mathbb B))^n$ is a row contraction such that its 
{\it Koszul complex}
\[
K^\cdot (M_z,H_m(\mathbb B))\stackrel{\epsilon_\lambda}{\longrightarrow} \mathbb C \rightarrow 0
\]
augmented by the point evaluations $\epsilon_\lambda:H_m(\mathbb B)\rightarrow \mathbb C,f\mapsto f(\lambda)$, at
arbitrary points $\lambda \in \mathbb B$, is exact \cite{GRS} (Proposition 2.6). In particular, the row operators
\[
H_m(\mathbb B)^n \longrightarrow H(\mathbb B),(h_i)^n_{i=1}\mapsto (\lambda-M_z)(h_i)^n_{i=1}=\sum^n_{i=1}(\lambda_i-M_{z_i})h_i
\]
have closed range and the operator-valued map $\mathbb B \rightarrow L(H_m(\mathbb B)^n,H(\mathbb B)),\lambda \mapsto \lambda-M_z$, is regular in the sense of \cite{M} (Theorem II.11.4).

The reciprocal of the kernel $K_m$ is given by the binomial sum
\[
K_m(z,w)^{-1}=\sum^m_{j=0}(-1)^j\binom{m}{j}\langle z,w\rangle^j.
\]
Since the row operator $M_z:H_m(\mathbb B)^n\rightarrow H_m(\mathbb B)$ has closed range, the operator $M^\ast_zM_z:{\rm Im}M^\ast_z\rightarrow 
{\rm Im}M^\ast_z$ is invertible. In \cite{E} it was shown that its inverse satisfies the identity
\[
(M^\ast_zM_z)^{-1} =\Big(\oplus \sum^{m-1}_{j=0}(-1)\binom{m}{j+1}\sigma^j_{M_z}(1_{H_m(\mathbb B)}\Big)|{\rm Im}M^\ast_z,
\]
where $\sigma_{M_z}(X)=\sum^n_{i=1}M_{z_i}XM^\ast_{z_i}$. We show that the commuting row contractions
 $T\in L(H)^n$ for which the operator-valued function $\mathbb B \rightarrow L(H^n,H)$, 
$\lambda  \mapsto  \lambda - T$, is regular at $z=0$ and which satisfy the operator identity 
\[
(T^\ast T)^{-1} = \Big(\oplus \sum^{m-1}_{j=0}(-1)^j\binom{m}{j+1}\sigma^j_T(1_H)\Big)|{\rm Im} T^\ast
\]
are precisely the commuting tuples which decompose into an orthogonal direct sum
\[
T=T_0\oplus T_1 \in L(H_0\oplus H_1)^n
\]
of a spherical coisometry $T_0\in L(H)^n$ and a tuple $T_1\in L(H_1)^n$ which is unitarily 
equivalent to the $m$-{\it shift} $M_z \in L(H_m(\mathbb B,\mathcal D))^n$ for some 
Hilbert space $\mathcal D$. 
We show that the coisometric part $T_0$ is absent if and only if
\[
\bigcap^\infty_{k=0} \sum_{|\alpha |=k}T^\alpha H=\lbrace 0 \rbrace.
\]
We thus extend corresponding one-variable results proved by Giselsson and Olafsson \cite{GO} for the standard weighted Bergman spaces on the unit disc to the case of the analytic Besov 
spaces $H_m(\mathbb B)$ on the unit ball.

For $m=1$, the space $H_1(\mathbb B)$ is the {\it Drury-Arveson space} and the validiy of the above operator identity means precisely that the row operator $T:H^n\rightarrow H$ is a partial 
isometry. Thus up to unitary equivalence, the commuting tuples $T \in L(H)^n$ that are regular at $z=0$
and for which the row operator $T: H^n \rightarrow H$ is a partial isometry, are precisely 
the direct sums $T=T_0 \oplus T_1 \in L(H_0\oplus H_1)^n$ of a spherical coisometry $T_0$ and a {\it Drury-Arveson shift} $T_1$. 
Specializing further to the case
$n=1$ one obtains Wold-type decompositions for partial isometries that contain the classical Wold decomposition theorem
and are closely related to corresponding results of Halmos and Wallen \cite{HW} for power partial isometries.

\vspace{1cm}

\centerline{\textbf{\S2 \, Analytic models}}

\vspace{.5cm}

Let $T \in L(H)^n$ be a commuting tuple of bounded operators on a complex Hilbert space $H$
such that $\sum_{1 \leq i \leq n} T_i H \subset H$ is a closed subspace. As usual we call
the space $W(T) = H \ominus \sum_{1 \leq i \leq n} T_i H_i$ the {\it wandering subspace} of $T$.
If the context is clear, we denote by $T$ also the induced row operator
$T: H^n \rightarrow H, (h_i)^n_{i=1} \mapsto \sum^n_{i=1} T_ih_i$, and
we write $T^*: H \rightarrow H^n, h \mapsto (T^*_i h)^n_{i=1}$, for its adjoint. 
Since $T: H^n\rightarrow H$ has closed range, the operator 
$T^* T: {\rm Im}\,T^* \rightarrow {\rm Im}\,T^*$ is invertible. We denote its inverse by
$(T^* T)^{-1}$.

Consider the column operator $L = (T^* T)^{-1}T^* \in L(H,H^n)$. 
Then $LT = P_{{\rm Im}\,T^*}$ and
\[
L(1_H - ZL)^{-1} (T- Z) = LT -  L \sum_{k=0}^{\infty} (ZL)^{k} Z(1_{H^n} - LT)
\]
\[
= P_{{\rm Im}\,T^*} - L(1_H - ZL)^{-1} Z P_{{\rm Ker}\,T}
\]
for $z \in \mathbb C^n$ with $\| z \| < 1/\| L \|$ $( = \infty$ if $L = 0)$. Here $Z: H^n \rightarrow H, 
(h_i)^n_{i=1} \mapsto \sum^n_{i=1} z_ih_i$, denotes the row operator induced by the 
complex $n$-tuple $z$ and $\| z \| = (\sum_{1 \leq i \leq n} |z_i |^2)^{1/2}$ is the Euclidean norm of $z$. 
Since ${\rm Im} \, L \subset {\rm Im} \, T^* 
= ({\rm Ker}\,T)^{\bot}$, it follows that $P(z) = (T-Z)L(1_H -ZL)^{-1} \in L(H)$ is a projection with
${\rm Im}\,P(z) \subset {\rm Im}(T - Z)$ for $z \in \mathbb C^n$ as above.
The identity $TL = T(T^*T)^{-1}T^* = P_{{\rm Im}\,T}$ yields that 
\[
1_H - P(z) = 1_H - (T-Z)L(1_H -ZL)^{-1} 
\]
\[
= (1_H - TL)(1_H -ZL)^{-1} = P_{W(T)} (1_H -ZL)^{-1} 
\]
and hence that $(1_H -P(z)) H = W(T)$ for $z \in \mathbb C^n$ with $\| z \| < 1/\| L \|$.

We call $T$ {\it regular} at $z = 0$ if there is a positive real number $\epsilon > 0$ 
such that, for $\| z \| < \epsilon$, the subspace
$(T - Z)H^n \subset H$ is closed and $H$ decomposes into the algebraic direct sum
\[
H = (T - Z)H^n \oplus W(T).
\]
Using this algebraic decomposition one obtains the identity
\[
(1_H - P(z))(T - Z) = 0 \quad {\rm for} \; \| z \| < \min(\epsilon, 1/ \| L \|).
\]
But then the identity theorem implies that ${\rm Im}(T-Z)  \subset {\rm Im}\, P(z)$ for
$\| z \| < 1/\| L \|$. Thus we find that
\[
 (1_H - P(z)) H = W(T) \quad  {\rm and} \quad P(z)H = (T - Z)H^n
\]
for $\| z \| < 1/\| L \|$.

The commuting tuple $T \in L(H)^n$ is regular at $z = 0$ if and only if the operator-valued
function $\mathbb C^n \rightarrow L(H^n,H)$, $z \mapsto T-Z$, is regular at $z = 0$ in the sense of \cite{M}
(see Definition II.10.21). This follows easily from the above remarks and Theorem II.11.4 in \cite{M}.
                             
Let us denote by $H^p(T,H) = {\rm Ker}\,\delta^p_T/{\rm Im}\,\delta^{p-1}_T$ the {\it cohomology groups}
of the Koszul complex $K^{\cdot}(T,H)$ of $T$ (see Section 2.2 in \cite{EP}).
It is well known and elementary to prove that the following conditions suffice to guarantee
the regularity of $T$ at $z = 0$.

\begin{lem} \label{regular}
Under either of the following three conditions:
\begin{enumerate}[{(i)}]
\item $T H^n = H$,
\item $T H^n \subset H$ is closed and $H^{n-1}(T,H) = \{0\}$,
\item there are an integer $N \geq 1$ and a real number $\delta > 0$ such that
      \[
			\dim H/(T-Z)H^n = N \quad {\rm for \; all} \; z \in B_{\delta}(0),
			\]
\end{enumerate}			
the tuple $T$ is regular at $z = 0$.
\end{lem}

\proof
We sketch the well known proofs. If $T H^n = H$, then $(T-Z)H^n = H$ for $z$ in a suitable neighbourhood 
of $z = 0$ and hence $T$ is regular at $z = 0$. Condition (ii) means precisely that the sequence
\[
\Lambda^{n-2}(\sigma,H) \xrightarrow{\binom{\delta^{n-2}_{z-T}}{\!\!\!\!\!\!\!\!\!0}}
\begin{matrix}
\Lambda^{n-1}(\sigma,H)\\[-0.2cm]
\oplus\\
W(T)
\end{matrix}
\xrightarrow{(\delta^{n-1}_{z-T},i)}\Lambda^n(\sigma,H) \longrightarrow 0,
\]
where $i: W(T) \hookrightarrow H$ denotes the inclusion map and the operators $\delta^p_{z-T}$
are the boundary maps of the Koszul complex of $z-T$ (Section 2.2 in \cite{EP}), is exact at $z = 0$. 
By Lemma 2.1.3 in \cite{EP} there is a positive real number $\epsilon < 1/\| L \|$ 
such that this sequence remains 
exact for every $z \in \mathbb C^n$ with $\| z \|< \epsilon$. But then
\[
(T-Z)H^n \oplus W(T) = H 
\]
and ${\rm Im}(T - Z) = P(z)H \subset H$ is closed for $\| z \| < \epsilon$.

Condition (iii) means that $T$ is a weak {\it Cowen-Douglas tuple} on $B_{\delta}(0)$ in the sense of
\cite{EL}. By the proof of Theorem 1.6 in \cite{EL} the tuple $T$ is regular at $z = 0$.
\proofend	

In the following let $T \in L(H)^n$ be a commuting tuple that is regular at $z = 0$.	We denote by	
$L_i \in L(H)$ $ (1 \leq i \leq n)$ the components of the column operator $L = (T^* T)^{-1}T^* \in L(H,H^n)$
and we use the notation
$L_i = L_{i_1} \cdots L_{i_k}$ for arbitrary index tuples $i = (i_1, \ldots ,i_k) \in \{1, \ldots ,n\}^k$.
To simplify the notation we write $\Omega_T = B_{1/\|L\|}(0)$ for the open Euclidean ball with 
radius $1/\|L\|$ at $z = 0$. We equip the space $\mathcal O(\Omega_T,W(T))$ of all analytic
$W(T)$-valued functions on $\Omega_T$ with its usual Fr\'{e}chet space topology of uniform convergence on
all compact subsets.

\begin{thm} \label{model}
Let $T \in L(H)^n$ be regular at $z = 0$.
Then the map
\[
V:H \rightarrow \mathcal O(\Omega_T,W(T)),\ (Vx)(z) = (1_H - P(z))x
\]
is continuous linear with $Vx \equiv x$ for $x \in W(T)$ and 
\begin{enumerate}[{(i)}]
\item $VT_i = M_{z_i}V$ \; $(i = 1, \ldots ,n)$,
\item ${\rm Ker}\,V = \bigcap^{\infty}_{k=0} \sum_{|\alpha|=k} T^{\alpha}H =\bigcap_{z \in \Omega_T} (T - Z)H^n$.
\end{enumerate}
\end{thm}

\proof 
By construction, for $z \in \Omega_T$ and $x \in H$, the vector 
\[
x(z) = (1_H -P(z))x  = P_{W(T)} (1_H - ZL)^{-1}x
\]
is the unique element  in $W(T)$ such that $x - x(z) \in {\rm Im}(T-Z).$
Obviously the vector $x(z)$ depends analytically on $z$ and the map V is continuous linear 
with $Vx \equiv x$ for $x \in W(T)$. Since for $z$ and $x$ as above,
\[
T_i x - z_ix(z) = T_i(x-x(z)) + (T_i-z_i)x(z) \in {\rm Im}(T-Z),
\]
the map $V$ intertwines the tuples $T$ on $H$ and $M_z$ on 
$\mathcal O(\Omega_T, W(T))$ componentwise. To calculate the kernel of $V$,
note that, for $x \in H$ and $z \in \Omega_T$,
\[
Vx(z) = \sum_{k=0}^{\infty}  P_{W(T)} (ZL)^k x
= \sum_{k=0}^{\infty} \sum_{| \alpha |=k} \Big( P_{W(T)} 
  \sum_{i \in I(\alpha)} L_i x \Big) z^{\alpha},
\]	
where for each $k \in \mathbb N$ and $\alpha \in \mathbb N^n$ with $|\alpha| = k$, the set $I(\alpha)$
consists of all index tuples $i = (i_1, \ldots ,i_k) \in \{1, \ldots ,n\}^k$ such that,
for each $j = 1, \ldots ,n$, exactly $\alpha_j$ of the indices $i_1, \ldots ,i_k$ 
equal $j$. The map $\Sigma_T: L(H) \rightarrow L(H), X \mapsto \sum^n_{i=1} T_i X L_i,$
is continuous linear with $P_{W(T)} = 1_H - TL = 1_H - \Sigma_T(1_H)$ and
\[
\sum_{j=0}^{k-1} \Sigma_T^j(P_{W(T)}) = 1_H - \Sigma_T^k(1_H) \quad (k \geq 0). 
\]
Hence for $x \in {\rm Ker}\,V$ and $k \geq 0$,
\[
0 = \sum_{j=0}^{k-1} \sum_{|\alpha|=j} T^{\alpha} \Big( P_{W(T)} \sum_{i \in I(\alpha)} L_i x\Big)
\]
\[ 
= \sum_{j=0}^{k-1} \Sigma_T^j(P_{W(T)})x = x - \sum_{|\alpha|=k} T^{\alpha} \Big(\sum_{i \in I(\alpha)} L_i x \Big).
\]
Thus ${\rm Ker}\,V \subset \bigcap^{\infty}_{k=0} \sum_{|\alpha|=k} T^{\alpha}H$. Conversely, if a vector
$x \in H$ belongs to the intersection on the right-hand side, then
\[
Vx \in \bigcap^{\infty}_{k=0} \sum_{|\alpha|=k} VT^{\alpha}H
\subset \bigcap^{\infty}_{k=0} \sum_{|\alpha|=k} M_z^{\alpha} \mathcal O(\Omega_T,W(T)) = \{ 0 \}.
\]
Thus the first equality in part (ii) has been shown. The second equality is obvious, since 
${\rm Ker}(1_H - P(z)) = {\rm Im}\,P(z) = (T - Z)H^n$ for all $z \in \Omega_T$.
\proofend

Elementary, even finite dimensional, examples show that Theorem \ref{model} need not be true
if instead of the regularity at $z = 0$ one only demands that the space $TH^n \subset H$ is closed.

Condition (ii) in Theorem \ref{model} implies that $W(T) \subset ({\rm Ker}\,V)^{\bot}$. An elementary
argument shows that $W(T)$ coincides with the wandering subspace 
of the compression of $T$ to $({\rm Ker}\,V)^{\bot}$.

In the following we use the notation $H_{\infty} = \bigcap^{\infty}_{k=0} \sum_{|\alpha|=k} T^{\alpha}H$. 
We call a commuting tuple $T \in L(H)^n$ {\it analytic} if $H_{\infty} = \{0\}$. 
If a commuting tuple $T \in L(H)^n$ is unitarily equivalent to the multiplication tuple
$M_z \in L(\mathcal H)^n$ on a functional Hilbert space $\mathcal H \subset \mathcal O(\Omega,\mathcal D)$
on a connected open zero neighbourhood $\Omega \subset \mathbb C^n$, then $T$ is necessarily  analytic. The next
result shows that, under the additional hypothesis that $T$ is regular at $z = 0$, also the converse implication holds.

Let $H/{\rm Ker}\, V \cong ({\rm Ker}\,V)^{\bot}$ be the quotient space of $H$ modulo the kernel of $V$. We denote 
the elements of $H/{\rm Ker}\, V$ by $x + {\rm Ker}\, V$.

\begin{cor} \label{analytic_model}
Let $T \in L(H)^n$ be regular at $z = 0$ and let
\[
V: H \rightarrow \mathcal O(\Omega_T,W(T))
\]
be the map from Theorem \ref{model}. Then $\mathcal H = {\rm Im}\,V \subset \mathcal O(\Omega_T,W(T))$ equipped
with the norm $\|x\| = \| x + {\rm Ker}\, V \|$ is a functional Hilbert space such that
\begin{enumerate}[{(i)}]
\item $P_{({\rm Ker}\,V)^{\bot}} T|({\rm Ker}\,V)^{\bot}$ is unitarily equivalent to $M_z \in L(\mathcal H)^n$
      via the unitary operator $V: ({\rm Ker}\,V)^{\bot} \rightarrow \mathcal H$,
\item the reproducing kernel $K_T: \Omega_T \times \Omega_T \rightarrow L(W(T))$ of $\mathcal H$ is given by
      \[
			K_T(z,w) = P_{W(T)}(1_H - ZL)^{-1}(1_H - L^*W^*)^{-1}|W(T).
			\]		
\end{enumerate}
\end{cor}

\proof
For $f \in \mathcal H$, there is a unique vector $x(f) \in ({\rm Ker}\,V)^{\bot}$ with $f = Vx(f)$. Since
$\lim_{k \rightarrow \infty} f_k = f$ in $\mathcal H$ if and only if $\lim_{k \rightarrow \infty} x(f_k) = x(f)$
in $H$, all point evaluations on $\mathcal H$ are continuous. Thus $\mathcal H$ is a functional Hilbert space.
For $y \in W(T)$ and $z \in {\rm Ker}\,V$,
\[
\langle (1_H - L^*W^*)^{-1}y,z \rangle = \langle y,(Vz)(w) \rangle = 0
\]
for every $w \in \Omega_T$. Let $f \in \mathcal H$, $y \in W(T)$ and $w \in \Omega_T$ be given. 
Define $x = x(f)$. Then 
\[
\langle f(w),y \rangle_{W(T)} = \langle P_{W(T)}(1_H - WL)^{-1}x,y \rangle_H = 
\langle x,(1_H - L^*W^*)^{-1}y \rangle_{({\rm Ker}\,V)^{\bot}}
\]
\[
= \langle Vx,V(1_H - L^*W^*)^{-1}y \rangle_{\mathcal H} = \langle f, K_T(\cdot,w)y \rangle_{\mathcal H}                     
\]
and hence $K_T$ is the reproducing kernel of the analytic functional Hilbert space $\mathcal H$. By construction
the compression of $T$ to $({\rm Ker}\,V)^{\bot}$ and $M_z \in L(\mathcal H)^n$ are unitarily equivalent 
via the unitary operator induced by $V$.
\proofend

In particular we obtain that each analytic tuple $T \in L(H)^n$ which is regular at $z = 0$ is unitarily
equivalent to a multiplication tuple $M_z \in L(\mathcal H)^n$ on a suitable analytic functional Hilbert space
$\mathcal H$ defined on a ball with center $0 \in \mathbb C^n$.
For single left invertible analytic operators, Corollary \ref{analytic_model} is due to Shimorin \cite{Sh}.

In the setting of Corollary \ref{analytic_model} the functional Hilbert space $\mathcal H \subset \mathcal O(\Omega_T,W(T))$
contains all polynomials $p(z) = \sum_{|\alpha| \leq N} x_{\alpha} z^{\alpha}$ with coefficients in $W(T)$. The
polnomials with coefficients in $W(T)$ are dense in $\mathcal H$ if and only if the space $({\rm Ker}\,V)^{\bot}$
is generated (as an invariant subspace) by the wandering subspace of $P_{({\rm Ker}\,V)^{\bot}} T|({\rm Ker}\,V)^{\bot}$.

\vspace{1cm}

\centerline{\textbf{\S3 \, Characterizations of Bergman shifts}}

\vspace{.5cm}

Let $T \in L(H)^n$ be a commuting tuple that is regular at $z = 0$.
As before we denote by $\sigma_T : \, L(H) \rightarrow L(H)$ the positive linear map
acting as $\sigma_T(X) = \sum_{1 \leq i \leq n} T_i X T^*_i$. 
We suppose in addition that $T$ satisfies the identity
\[
(T^* T)^{-1} = (\oplus \Delta_T)|{\rm Im}\,T^*,
\]
where $\Delta_T \in L(H)$ is the operator defined by
\[
\Delta_T = \sum^{m-1}_{j=0}(-1)^j\binom{m}{j+1}\sigma^j_T(1_H).
\]
Let us define $\delta_T \in L(H)$ by
\[
\delta_T = ({\rm Im}\,T\stackrel{T^*}{\longrightarrow}{\rm Im}\,T^*)^{-1}(T^* T)^{-1}T^*.
\]
Then ${\rm Im}\,\delta_T = {\rm Im}\,T$ and
\[
T^*_i \delta_T = \Delta_T T^*_i \quad \quad (i=1,\ldots,n).
\]
Using these intertwining relations, we find that $(\Delta_T T^*_i)(\Delta_T T^*_j) = 
\Delta_T T^*_i T^*_j \delta_T = (\Delta_T T^*_j)(\Delta_T T^*_i)$ for $i,j = 1, \ldots, n$. In the
following we use the same notation for the column operator $L: H \rightarrow H^n, 
x \mapsto (T^* T)^{-1}T^* x = (\Delta_T T^*_i x)^n_{i=1},$ and the commuting tuple 
$L = (\Delta_TT^*_i)^n_{i=1} \in L(H)^n$. 
Since $L$ is commuting, the representation of the map $V$ obtained in the proof of Theorem \ref{model}
simplifies to
\[
(Vx)(z) = \sum_{\alpha \in \mathbb N^n} \gamma_{\alpha} (P_{W(T)}L^{\alpha}x) z^{\alpha} \quad (x \in H, \; z \in \Omega_T),
\]
where $\gamma_{\alpha} = |\alpha|!/\alpha !$ for $\alpha \in \mathbb N^n$.
 
\begin{lem} \label{cancel}
For $\alpha, \beta \in \mathbb N^n$, we have
\[
\gamma_{\alpha} P_{W(T)} L^{\alpha} T^{\beta} = \gamma_{\alpha - \beta} P_{W(T)} L^{\alpha - \beta},
\]
where the right-hand side has to be read as zero whenever $\alpha - \beta$ has negative components.
\end{lem}

\proof
For $\beta \in \mathbb N^n,x \in H$ and $ z \in \Omega_T$,
\[
\sum_{\alpha \in \mathbb N^n} \gamma_{\alpha} P_{W(T)} (L^{\alpha}x) z^{\alpha + \beta}
= z^{\beta} (Vx)(z) = (VT^{\beta}x)(z) = \sum_{\alpha \in \mathbb N^n} \gamma_{\alpha} P_{W(T)} 
(L^{\alpha} T^{\beta}x) z^{\alpha}.
\]
The proof follows by comparing the coefficients of these convergent power series.
\proofend

Let us apply the above constructions to the particular case of the multiplication tuple $T = M_z \in L(H_m(\mathbb B))^n$.
By definition $H_m(\mathbb B)$ is the analytic functional Hilbert space  with
reproducing kernel $K_m:\mathbb B \times \mathbb B \rightarrow \mathbb C,\ K_m(z,w) 
= (1 - \langle z,w \rangle)^{-m}$. We consider only the case where the exponent $m \geq 1$ is a positive integer. 
The commuting tuple
$M_z = (M_{z_1},\ldots,M_{z_n})\in L(H_m(\mathbb B))^n$, consisting of the multiplication operators 
$M_{z_i}: H_m(\mathbb B) \rightarrow H_m(\mathbb B),\ f \mapsto z_i f$, with the coordinate functions, is regular at $z = 0$,
Indeed, for each point $\lambda \in \mathbb B$, the Koszul complex
$K^{\cdot}(\lambda - M_z,H_m(\mathbb B))$ is exact in degree $p = 0, \ldots, n-1$ and
$\dim H^n(K^{\cdot}(\lambda - M_z,H_m(\mathbb B))) = 1$ (see e.g. Proposition 2.6 in \cite{GRS}).

By Lemma 1 and Lemma 3 in \cite{E} (see also the proof of Lemma 3 in \cite{E}) we know that
\[
(M^*_z M_z)^{-1} = (\oplus \Delta_{M_z})|{\rm Im}\,M^*_z,
\]
where $\Delta_{M_z} = \sum_{j=0}^{m-1} (-1)^j \binom{m}{j+1} \sigma^j_{M_z}(1_{H_m(\mathbb B)})$ acts as the
diagonal operator
\[
\Delta_{M_z} \sum_{k=0}^{\infty} (\sum_{|\alpha| = k} f_{\alpha} z^{\alpha}) = 
\sum_{k=0}^{\infty} \frac{m+k}{1+k} (\sum_{|\alpha| = k} f_{\alpha} z^{\alpha}).
\]
In this particular example, $W(M_z) = \mathbb C$ and $P_{W(M_z)} \in  L(H_m(\mathbb B))$ is the
orthogonal projection onto the closed subspace $\mathbb C \subset  H_m(\mathbb B)$ consisting of all constant
functions. Furthermore, the intertwining relation $M^*_z \delta = (\oplus \Delta_{M_z}) M^*_z$ holds with the
diagonal operator $\delta = \delta_{M_z}: H_m(\mathbb B) \rightarrow H_m(\mathbb B),$
\[
\delta \sum_{k=0}^{\infty} (\sum_{|\alpha| = k} f_{\alpha} z^{\alpha}) = f_0 +
\sum_{k=1}^{\infty} \frac{m+k-1}{k} (\sum_{|\alpha| = k} f_{\alpha} z^{\alpha})
\]
(see the proof of Lemma 3 in \cite{E}).

\begin{lem} \label{concrete_inter}
The commuting tuple $L_{M_z} = (\Delta_{M_z}M^*_{z_1}, \ldots , \Delta_{M_z}M^*_{z_n}) \in L(H_m(\mathbb B))^n$
satisfies the identities
\[
P_{W(M_z)} L^{\alpha}_{M_z} = \binom{m + |\alpha| - 1}{|\alpha|} P_{W(M_z)} M_z^{* \alpha} \quad (\alpha \in \mathbb N^n).
\]
\end{lem}

\proof
For $\alpha = 0$, the identity obviously holds. Suppose that the result has been shown for each multiindex
$\alpha \in \mathbb N^n$ with $|\alpha| \leq k$ and fix an $\alpha \in \mathbb N^n$ with $|\alpha| = k$
as well as an index $ i \in \{1, \ldots ,n\}$. Then
\begin{align*}
P_{W(M_z)} L^{\alpha + e_i}_{M_z} &= \binom{m + k - 1}{k} P_{W(M_z)} M_z^{* \alpha} \Delta_{M_z} M^*_{z_i}\\
&= \binom{m + k - 1}{k} P_{W(M_z)} M_z^{* (\alpha +e_i)} \delta\\
&= \binom{m + k - 1}{k} \frac{m+k}{k+1} P_{W(M_z)} M_z^{* (\alpha +e_i)}\\
&= \binom{m + k}{k + 1}  P_{W(M_z)} M_z^{* (\alpha + e_i)}.
\end{align*}
Thus the assertion follows by induction on $|\alpha|$.
\proofend

We use the result proved in Lemma \ref{concrete_inter} for $M_z \in L(H_m(\mathbb B))^n$ to prove the correspon-ding
result for the commuting tuple $T$ fixed at the beginning of Section 3. 

\begin{lem}\label{intertwining}
For $\alpha \in  \mathbb N^n$, the identity
\[
P_{W(T)} L^{\alpha} = \binom{m + |\alpha| - 1}{|\alpha|} P_{W(T)} T^{* \alpha}
\]
holds.
\end{lem}

\proof
Again we use induction on $|\alpha|$. Suppose that the result holds for $|\alpha| \leq k$.
Let $\alpha \in \mathbb N^n$ be a multiindex with  $|\alpha| = k$ and let 
$i \in \{1, \ldots ,n \}$ be arbitrary. Using Lemma \ref{cancel} and the induction hypothesis we obtain
\begin{align*}
P_{W(T)} L^{\alpha + e_i} &= P_{W(T)} L^{\alpha} \Delta_T T_i^*\\
&= P_{W(T)} L^{\alpha} \sum_{j=0}^{m-1} (-1)^j \binom{m}{j+1} \sum_{|\beta|=j} \gamma_{\beta} T^{\beta} T^{* \beta} T^*_i\\
&= \sum_{j=0}^{m-1} (-1)^j \binom{m}{j+1} \sum_{\substack{|\beta|= j\\ \alpha \geq \beta}} 
   \frac{\gamma_{\beta} \gamma_{\alpha - \beta}}{\gamma_{\alpha}}
   P_{W(T)} L^{\alpha - \beta} T^{* (\beta + e_i)}\\
&= \Big( \sum_{j=0}^{m-1} (-1)^j \binom{m}{j+1} \sum_{\substack{|\beta|= j\\ \alpha \geq \beta}} 
   \frac{\gamma_{\beta} \gamma_{\alpha - \beta}}{\gamma_{\alpha}}
   \binom{m + |\alpha - \beta| -1}{|\alpha - \beta|}\Big) P_{W(T)} T^{* (\alpha + e_i)}.	 
\end{align*}

Here by definition $\alpha \geq \beta$ means that $\alpha_i \geq \beta_i$ for $i = 1, \ldots ,n$.
Next observe that the preceding chain of equalities remains true if $T, \Delta_T$ and $L$ are
replaced by $M_z, \Delta_{M_z}$ and $L_{M_z}$. But in this case we know from Lemma \ref{concrete_inter}
that
\[
P_{W(M_z)} L^{\alpha + e_i}_{M_z} = \binom{m + |\alpha|}{|\alpha| + 1} P_{W(M_z)} M^{* (\alpha + e_i)}_z.
\]
Hence we may conclude that
\[
\sum_{j=0}^{m-1} (-1)^j \binom{m}{j+1} \sum_{\substack{|\beta|= j\\ \alpha \geq \beta}} 
\frac{\gamma_{\beta} \gamma_{\alpha - \beta}}{\gamma_{\alpha}}
\binom{m + |\alpha - \beta| -1}{|\alpha - \beta|} = \binom{m + |\alpha|}{|\alpha| + 1}.
\]
This observation completes the inductive proof.
\proofend

Since $V: H \rightarrow \mathcal O(\Omega_T,W(T)), Vx(z) = 
\sum_{\alpha \in \mathbb N^n} \gamma_{\alpha} (P_{W(T)} L^{\alpha} x) z^{\alpha},$
is a continuous linear map that intertwines the tuples $T$ on $H$ and $M_z$ on
$\mathcal O(\Omega_T,W(T))$ componentwise, the kernel of $V$ is a closed invariant subspace
for $T$. Much more than this is true.

\begin{lem}\label{reducing}
The kernel of $V$ is reducing for $T$ with
\begin{enumerate}
\item ${\rm Ker}\,V = H_{\infty}
      = \{ x \in H; P_{W(T)} T^{* \alpha} x = 0 \; {\rm for \, all} \; \alpha \in \mathbb N^n \}$,
\item	$({\rm Ker}\,V)^{\bot} = \bigvee_{\alpha \in \mathbb N^n} T^{\alpha}W(T)$.		
\end{enumerate}
\end{lem}

\proof
The first equality in part (a) holds by Theorem \ref{model}. Since
\[
Vx(z) = \sum_{\alpha \in \mathbb N^n} \gamma_{\alpha} (P_{W(T)} L^{\alpha} x) z^{\alpha}
\]
for $x \in H$ and $ z \in \Omega_T$, the kernel of $V$ consists of all vectors $x \in H$ with
$P_{W(T)} L^{\alpha} x = 0$ for all $\alpha \in \mathbb N^n$. Thus the second equality in part (a)
follows from Lemma \ref{intertwining}. In view of the identity
\[
\langle x, T^{\alpha} y \rangle = \langle P_{W(T)} T^{* \alpha} x, y \rangle \quad (x \in H, y \in W(T), \alpha \in \mathbb N^n)
\]
part (b) follows from (a). Both parts together imply that ${\rm Ker}\,V$ is a reducing subspace for $T$.
\proofend

In the following we write $[M] \subset H$ for the smallest closed linear subspace of $H$ which
contains a given subset $M \subset H$.
For a complex Hilbert space $\mathcal E$, we denote by $H_m(\mathbb B,\mathcal E)$ the $\mathcal E$-valued 
analytic functional Hilbert space with reproducing kernel
\[
K^\mathcal E_m:\mathbb B \times \mathbb B \rightarrow L(\mathcal E), K^{\mathcal E}_m(z,w) = \frac{1_\mathcal E}{(1-\langle z,w\rangle)^m}
\]
on $\mathbb B$. A well known alternative description of the space $H_m(\mathbb B,\mathcal E)$ is given by
\[
H_m(\mathbb B,\mathcal E) = \{ f = \sum_{\alpha \in \mathbb N^n} f_\alpha z^\alpha \in \mathcal O(\mathbb B,\mathcal E); \| f \|^2
= \sum_{\alpha \in \mathbb N^n}\frac{\| f_\alpha \|^2}{\rho_m(\alpha)} < \infty \},
\]
where
$
\rho_m(\alpha) = \frac{(m + |\alpha|-1)!}{\alpha!(m - 1)!}.
$

\begin{thm} \label{shift}
Let $T \in L(H)^n$ be a commuting tuple that is regular at $z = 0$ and satisfies the identity
$(T^* T)^{-1} = (\oplus \Delta_T)|{\rm Im}\,T^*$.
Then the map
\[
U: [W(T)] \rightarrow H_m(\mathbb B,W(T)), \; Ux(z) = \sum_{\alpha \in \mathbb N^n} \gamma_{\alpha} (P_{W(T)} L^{\alpha} x) z^{\alpha}
\]
is a unitary operator which componentwise intertwines the tuples $T|[W(T)]$ and $M_z \in$  $L(H_m(\mathbb B,W(T)))^n$.
\end{thm}

\proof
For $N \in \mathbb N$ and $x_{\alpha} \in W(T)$ $(|\alpha| \leq N)$, we have
\[
\| \sum_{|\alpha| \leq N} T^{\alpha} x_{\alpha} \|^2 = \sum_{|\alpha|,|\beta| \leq N} 
\langle P_{W(T)} T^{ *\beta}T^{\alpha} x_{\alpha}, x_{\beta} \rangle
= \sum_{|\alpha|,|\beta| \leq N} \langle x_{\alpha}, P_{W(T)} T^{ *\alpha}T^{\beta}  x_{\beta} \rangle.
\]
Using first Lemma \ref{intertwining} and then twice Lemma \ref{cancel}, we find that
\begin{align*}
\| \sum_{|\alpha| \leq N} T^{\alpha} x_{\alpha} \|^2 &= 
\sum_{|\alpha| \leq N} \binom{m + |\alpha| -1}{|\alpha|}^{-1} \langle P_{W(T)} L^{\alpha}T^{\alpha} x_{\alpha}, x_{\alpha} \rangle\\
&= \sum_{|\alpha| \leq N} \Big( \binom{m + |\alpha| -1}{|\alpha|} \gamma_{\alpha} \Big)^{-1} \| x_{\alpha} \|^2\\
&= \sum_{|\alpha| \leq N} \frac{\| x_{\alpha} \|^2}{\rho_m(\alpha)} = 
   \| \sum_{|\alpha| \leq N} x_{\alpha} z^{\alpha} \|^2_{H_m(\mathbb B,W(T))}.
\end{align*}
Since the polynomials with coefficients in $W(T)$ are dense in $H_m(\mathbb B,W(T))$,
there is a unique unitary operator $U: [W(T)] \rightarrow H_m(\mathbb B,W(T))$ with
$U(\sum_{|\alpha| \leq N} T^{\alpha} x_{\alpha}) = \sum_{|\alpha| \leq N} x_{\alpha} z^{\alpha}$ for all finite
families $(x_{\alpha})_{|\alpha| \leq N}$ in $W(T)$. In particular it follows that, for  
$h \in {\rm span}\{T^{\alpha}x; \; \alpha \in \mathbb N^n \;{\rm and} \; x \in W(T) \}$, 
the analytic functions $Uh \in \mathcal O(\mathbb B,W(T))$ and $Vh \in \mathcal O(\Omega_T,W(T))$
have the same Taylor coefficients at $z = 0$. The continuity of the maps
$U: [W(T)] \rightarrow \mathcal O(\mathbb B,W(T))$ and $V: [W(T)] \rightarrow \mathcal O(\Omega_T,W(T))$
implies that $Ux$ and $Vx$ have the same Taylor coefficients at $z = 0$ for every $x \in [W(T)]$.
But then
\[
Ux(z) = \sum_{\alpha \in \mathbb N^n} \gamma_{\alpha} (P_{W(T)} L^{\alpha} x) z^{\alpha}
\]
for $x \in [W(T)]$ and $z \in \mathbb B$. Since $V$ intertwines $T \in L(H)^n$ and $M_z$ on $\mathcal O(\Omega_T,W(T))$,
the identity theorem implies that $U$ satisfies the same intertwining relation. 
\proofend

For a commuting tuple $S \in L(H)^n$ its $k$th {\it order defect operators} are defined by
\[
\Delta_S^{(k)} = (I - \sigma_S)^k(1_H) = \sum^k_{j=0} (-1)^j \binom{k}{j} \sigma_S^j(1_H) \quad (k \in \mathbb N).
\]
The tuple $S$ is called an $m$-{\it hypercontraction} if $\Delta_S^{(1)} \geq 0$ and $\Delta_S^{(m)} \geq 0$. A commuting
tuple $S \in L(H)^n$ is said to be of {\it type} $C_{\cdot 0}$ if SOT-$\lim_{k \rightarrow \infty} \sigma_S^k(1_H) = 0$.

\begin{cor} \label{pure} 
Let $T \in L(H)^n$ be as in Theorem \ref{shift}. The following conditions on $T$ are equivalent:
\begin{enumerate}[{(i)}]
\item $T$ is analytic,
\item $\|x\|^2 = \| \sum_{\alpha \in \mathbb N^n} \gamma_{\alpha} (P_{W(T)} L^{\alpha} x) z^{\alpha} \|^2_{H_m(\mathbb B,W(T))}$
      for all $x \in H$, 
\item $T$ is of type $C_{\cdot 0}$,			
\item $T$ is unitarily equivalent to $M_z \in L(H_m(\mathbb B,\mathcal D))^n$ for some Hilbert space $\mathcal D$.		
\end{enumerate}
\end{cor}

\proof
The equivalence of (i) and (ii) follows from Lemma \ref{reducing} and Theorem \ref{shift}. The implication
(i) to (iv) follows from Theorem \ref{shift}. It is well known that $M_z \in L(H_m(\mathbb B,\mathcal D))^n$
satisfies the $C_{\cdot 0}$-condition
\[
{\rm SOT}-\lim_{k \rightarrow \infty} \sigma^k_{M_z}(1_{H_m(\mathbb B,\mathcal D)}) = 0.
\]
Since this condition is preserved under unitary equivalence, the implication (iv) to (iii) holds.

Let us suppose that $T$ satisfies condition (iii). To complete the proof note first that
\begin{align*}
P_{W(T)} &= 1_H - T(T^*T)^{-1} T^* = 1_H - T(\oplus \Delta_T) T^*\\
&= 1_H - \sigma_T \Big( \sum_{k=0}^{m-1} (-1)^k \binom{m}{k+1} \sigma_T^k(1_H)\Big)\\
&= 1_H + \sum_{k=1}^m (-1)^k \binom{m}{k} \sigma_T^k(1_H)\\ 
&= (I - \sigma_T)^m(1_H).
\end{align*}
It is well known that a commuting tuple $T \in L(H)^n$ of type $C_{\cdot 0}$ for which the $m$th order 
defect operator $\Delta^{(m)}_T = (I - \sigma_T)^m(1_H)$ 
is positive is an $m$-hypercontraction
(see \cite{E, MV}). Thus it follows from the dilation theory for  $m$-hypercontractions
(see e.g. \cite{E}) that the map
\[
j: H \rightarrow H_m(\mathbb B,H) \;
jx = \sum_{\alpha \in \mathbb N^n} \rho_m(\alpha) \big((\Delta^{(m)}_T)^{1/2} T^{* \alpha} x\big) z^{\alpha}
\]
defines an isometric intertwiner between the tuples $T^*$ on $H$ and $M^*_z$ on $H_m(\mathbb B,H)$.
Using Lemma \ref{intertwining} we find that
\[
(Vx)(z) = \sum_{\alpha \in \mathbb N^n} \gamma_{\alpha} (P_{W(T)} L^{\alpha} x) z^{\alpha}
        = \sum_{\alpha \in \mathbb N^n} \rho_m(\alpha) (P_{W(T)} T^{* \alpha} x) z^{\alpha} = (jx)(z)
\]
for $x \in H$ and $z \in \Omega_T$. Hence ${\rm Ker}\,V = {\rm Ker}\,j = \{0\}$ and the proof is complete.
\proofend

Let $S\in L(H)^n$ be an $m$-hypercontraction. Since 
\[
0\leq \sigma^{k+1}_S(1_H)\leq \sigma^k_S(1_H)\leq 1_H\quad \quad (k\geq 0),
\]
the strong limit $S_\infty = {\rm SOT}-\lim_{k\rightarrow \infty}\sigma^k_S(1_H)$ exists. 
It is well known that the map $j:H\rightarrow H_m(\mathbb B,H)$,
\[
jx=\sum_{\alpha \in \mathbb N^n}\rho_m(\alpha)((\Delta^{(m)}_S)^{1/2}S^{\ast \alpha}x)z^\alpha
\]
defines a contraction that intertwines the tuples $S^\ast \in L(H)^n$ and $M^\ast_z \in L(H_m(\mathbb B,H))^n$. More precisely one can show that
\[
\| jx\|^2+\| S_\infty x\|^2=\| x \|^2
\]
for every vector $x \in H$.

Let us recall that a commuting tuple $S \in L(H)^n$ is called a {\it spherical isometry} if $\sum_{1 \leq i \leq n} S^*_i S_i = 1_H$
or, equivalently, if $\sum_{1 \leq i \leq n} \| S_i x \|^2 = \| x \|^2$ for each vector $x \in H$. By a result of Athavale
\cite{Ath} each spherical isometry $S \in L(H)^n$ is subnormal and its minimal normal extension is a {\it spherical unitary},
that is, a commuting tuple $N \in L(K)^n$ of normal operators such that $\sum_{1 \leq i \leq n} N^*_i N_i = 1_K$. A
{\it spherical coisometry} is a commuting tuple $S \in L(H)^n$ such that its adjoint $S^* \in L(H)^n$ is a spherical
isometry.

\begin{thm}\label{row}
Let $T \in L(H)^n$ be a commuting row contraction that is regular at $z=0$. Then $T$ satisfies the operator identity
\[
(T^* T)^{-1} = \Big(\oplus \sum^{m-1}_{j=0}(-1)^j\binom{m}{j+1}\sigma^j_T(1_H)\Big)|{\rm Im}\, T^*
\]
if and only if $T = T_0 \oplus T_1 \in L(H_0\oplus H_1)^n$ is the direct sum of a spherical coiso-metry
$T_0 \in L(H_0)^n$ and a tuple $T_1 \in L(H_1)^n$ which
is unitarily equivalent to $M_z \in L(H_m(\mathbb B,\mathcal D))^n$ 
for some Hilbert space $\mathcal D$.
\end{thm}

\proof
Suppose that $T$ satisfies the above operator identity. Then by Lemma \ref{reducing} the space $H$ is
the orthogonal sum $H = H_{\infty} \oplus [W(T)]$ of closed subspaces reducing $T$. According to 
Theorem \ref{shift} the restriction $T_1 = T|[W(T)]$ is unitarily equivalent to $M_z \in L(H_m(\mathbb B,W(T)))^n$.
The proof of Corollary \ref{pure} shows that
\[
\Delta_T^{(m)} = (1 - \sigma_T)^m(1_H) = P_{W(T)} \geq 0.
\]
Since $T$ is a row contraction, also $\Delta_T^{(1)} = 1_H - \sum_{1 \leq i \leq n} T_i T^*_i \geq 0$.
Thus the tuple $T$ is an $m$-hypercontraction. Using the remarks following Corollary \ref{pure} as well as Lemma \ref{intertwining}
we obtain that the map $j: \, H \rightarrow H_m(\mathbb B,H)$,
\[
jx = \sum_{\alpha \in \mathbb N^n} \rho_m(\alpha) (P_{W(T)} T^{* \alpha} x) z^{\alpha}
    =\sum_{\alpha \in \mathbb N^n} \gamma_{\alpha} (P_{W(T)} L^{\alpha} x) z^{\alpha}
\]
is a well defined contraction with $\| jx \|^2 + \langle T_{\infty}x,x \rangle = \| x \|^2$ for all $x \in H$.
Since $T_{\infty} \leq \sigma_T(1_H) \leq 1_H$, it follows that
\[
\sum_{i=1}^n \| T^*_i x \|^2 = \langle \sigma_T(1_H)x,x \rangle = \| x \|^2
\]
for $x \in {\rm Ker}\,j = {\rm Ker}\,V = H_{\infty}$. Thus $(T|H_{\infty})^*$ is a spherical isometry.

To prove the converse let us first consider a commuting tuple $T \in L(H)^n$ such that $T^*$ is a spherical
isometry. Since $\sigma_T(1_H) = TT^* =1_H$, it follows that $T^* T = P_{{\rm Im}\,T^*}$ and
\[
\sum^{m-1}_{j=0}(-1)^j\binom{m}{j+1}\sigma^j_T(1_H) = \Big(\sum^{m-1}_{j=0}(-1)^j\binom{m}{j+1}\Big) 1_H = 1_H.
\]
On the other hand, by Lemma 1 and Lemma 3 in \cite{E} it follows that, for any Hilbert space $\mathcal D$, 
the tuple $M_z \in L(H_m(\mathbb B,\mathcal D))^n$ satisfies the operator identity
\[
(M^*_zM_z)^{-1} = \big(\oplus \sum^{m-1}_{j=0}(-1)^j\binom{m}{j+1}\sigma^j_{M_z}(1_{H_m(\mathbb B,\mathcal D)})\big)|{\rm Im}M^*_z.
\]
Since the validity of this operator identity is preserved under unitary equivalence and the passage to direct sums,
also the reverse implication follows.
\proofend

In the single-variable case $n = 1$ Theorem \ref{row} implies that the left invertible contractions
that satisfy the operator identity
\[
(T^* T)^{-1} = \sum^{m-1}_{j=0}(-1)^j\binom{m}{j+1}\sigma^j_T(1_H)
\]
are precisely the operators $T \in L(H)$ that decompose into the orthogonal direct sum $T = T_0 \oplus T_1$ of
a unitary operator $T_0$ and  an operator $T_1$  which is unitarily equivalent to an $m$-shift $M_z \in L(H_m(\mathbb D,\mathcal D))$.
This is a slight variant of the main result of \cite{GO}.

In the particular case $m = 1$ the result stated in Theorem \ref{row} takes the form.

\begin{cor} \label{m=1}
Let $T \in L(H)^n$ be a commuting tuple that is regular at $z = 0$. Then $T: H^n \rightarrow H$ is a partial
isometry if and only if $T = T_0 \oplus T_1 \in L(H_0\oplus H_1)^n$ is the direct sum of a spherical coisometry
$T_0 \in L(H_0)^n$ and a tuple  $T_1 \in L(H_1)^n$ which is 
unitarily equivalent to $M_z \in L(H_1(\mathbb B,\mathcal D))^n$ 
for some Hilbert space $\mathcal D$.
\end{cor}

\proof 
For $m = 1$, we have $ \sum^{m-1}_{j=0}(-1)^j\binom{m}{j+1}\sigma^j_T(1_H) = 1_H$.
Thus in this case the operator identity from Theorem \ref{row} means precisely that 
$T: H^n \rightarrow H$ is a partial isometry. Hence the assertion follows immediately from Theorem \ref{row}.
\proofend

In the case $m = 1 = n$ the preceding results yield a Wold decomposition for partial isometries $T \in L(H)$
that are regular at $z = 0$ which contains the classical Wold decomposition for isometries. Corollary \ref{m=1}
implies that all powers $T^k$ of partial isometries $T \in L(H)$ that are regular at $z = 0$ are partial isometries
again. Thus partial isometries that are regular at $z = 0$ are {\it power partial isometries} in the sense of 
Halmos and Wallen \cite{HW}. In \cite{HW} a Wold-decomposition theorem for general power partial isometries is proved.

Since a non-unitary partial isometry on a finite dimensional Hilbert space cannot admit a decomposition as 
in Corollary \ref{m=1}, Corollary \ref{m=1} and Theorem \ref{row} cannot be expected to hold without the
hypothesis that the given tuple $T$ is regular at $z = 0$.

\begin{minipage}[t]{0.5\textwidth}
J. Eschmeier \\
Fachrichtung Mathematik\\
Universit\"at des Saarlandes\\
Postfach 151150\\
D-66041 Saarbr\"ucken, Germany
\end{minipage}

\vspace{.4cm}

\begin{minipage}[t]{0.5\textwidth}
S. Langend\"orfer\\
Fachrichtung Mathematik\\
Universit\"at des Saarlandes\\
Postfach 151150\\
D-66041 Saarbr\"ucken, Germany
\end{minipage}

\texttt{eschmei@math.uni-sb.de, langendoerfer@math.uni-sb.de}

\end{document}